\documentclass[11pt]{article}
\usepackage[centertags]{amsmath}
\usepackage{amsmath}
\usepackage{hyperref}
\usepackage{amssymb}
\usepackage{amsthm}
\usepackage{latexsym}
\usepackage{times}
\usepackage[T1]{fontenc}
\usepackage[latin1]{inputenc}
\usepackage{color}
\usepackage{fancyhdr}
\usepackage{cite}
\usepackage{indentfirst}
\DeclareMathOperator*{\esssup}{ess\,sup}

\numberwithin{equation}{section}

\newcommand{\be}{\begin{equation}}
\newcommand{\ee}{\end{equation}}
\def \ed {\end{document}}

\newtheorem{theorem}{Theorem}[section]
\newtheorem{proposition}{Proposition}[section]
\newtheorem{lemma}{Lemma}[section]
\newtheorem{remark}{Remark}[section]

\def \esssup {\mbox{ess sup}}

\def \R{\mathbb{R}}
\def \P{\mathbb{P}}

\def \E{\mathbb{E}}

\def \F{\mathcal{F}}

\def \I{\mathcal{I}}
\def \A{\mathcal{A}}
\def \L{\mathcal{L}}

\def \D{\mathcal{D}}
\def \Y{\mathcal{Y}}
\def \p{\mathcal{P}}
\def \bf{\textbf}

\begin{document}

\title{Infinite Horizon Multi-Dimensional BSDE with Oblique Reflection and Switching Problem}


\author{ Brahim EL ASRI \thanks{Universit\'e Ibn Zohr, Equipe. Aide \`a la decision,
ENSA, B.P.  1136, Agadir, Maroc. e-mail: b.elasri@uiz.ac.ma } \,\,\, and \, Nacer OURKIYA \thanks{Universit\'e Ibn Zohr, Equipe. Aide \`a la decision,
ENSA, B.P.  1136, Agadir, Maroc. e-mail: nacer.ourkiya@edu.uiz.ac.ma.}}
\maketitle

\begin{abstract}
This paper studies a system of multi-dimensional reflected backward stochastic differential equations with oblique reflections (RBSDEs for short) in infinite horizon associated  to switching problems.
 The existence and uniqueness of the adapted solution is obtained by using a method based on a combination of penalization, verification and contraction property.
\end{abstract}

\\


\medskip
\section{Introduction}
In this paper we study a system of multi-dimensional RBSDE with oblique reflection in infinite horizon.
\par For $i\in\mathcal{I}:=\{1,...,m\}$ and $t\geq0$, we define the multi-dimensional RBSDE with oblique reflections by, $\forall r\in\R^+$
\begin{equation}\label{rbsde}
\begin{cases}
e^{-rt}Y^{i}_t = \int_{t}^{+\infty} e^{-rs}f_i(s,X_s,Y^{1}_s,...,Y^{m}_s,Z^{i}_s) \, ds + K^{i}_{\infty} - K^{i}_t-  \int_{t}^{+\infty}e^{-rs}Z^{i}_s \, dB_s;\\
\lim_{t \to+\infty} e^{-rt}Y^{i}_t=0,\\
\forall t \geq 0, \quad e^{-rt}Y^{i}_t \geq e^{-rt} \underset{j\in\mathcal{I}^{-i}}{\max}(Y^{j}_t-g_{ij}),\\
\int_{0}^{+\infty} e^{-rs}\{Y^{i}_s - \underset{j\in\mathcal{I}^{-i}}{\max}(Y^{j}_s-g_{ij})\} \, dK^i_s = 0,
\end{cases}
\end{equation}
where $\mathcal{I}^{-i}:=\mathcal{I}\textbackslash\{i\}$, $B$ is a standard Brownian motion on a complete probability space $(\Omega,\mathcal{F},\P)$, $f_i$ the generator which is Lipschitz continuous w.r.t. $y$ and $z$, $X:=(X_s)_{s\geq0}$ be a $\p$-measurable, $\R^k$-valued continuous stochastic process and $g_{ij}$ represent the switching costs from mode $i$ to mode $j.$ We aim at finding a $m$-triples of $(\mathcal{F}_t)_{t\geq0}$-adapted processes $(Y^i,Z^i,K^i)_{i\in\I}$, which solves $\P$-a.s. the system of multi-dimensional RBSDE $\eqref{rbsde}$.
 \par One-dimensional RBSDEs were first studied by El Karoui et al. \cite{ppq} in finite horizon case. Later, Hamadène et al. \cite{HLW} proved existence and uniqueness results of the solution for infinite horizon RBSDEs. The authors applied these results
to get the existence of optimal control strategy for the mixed control problem. The literature on this specific form of equation has then grown very importantly due to their range of application, in particular  in the field of stochastic control or mathematical finance.
\par Multi-dimensional RBSDEs were studied by Gegout-Petit and Pardoux \cite{GP}, but their BSDE
is reflected on the boundary of a convex domain along the inward normal direction, and their method depends heavily on the properties of this inward normal reflection (see $(1)$-$(3)$ in \cite{GP}). And then, by Ramasubramanian \cite{Ram} in an orthant with some restriction on the direction of oblique reflection and the driver $f$.
\par Another type of multi-dimensional RBSDEs occurs in the context of the optimal switching problem (see e.g. \cite{BCR,CEK,CEK2,CR,HJ,HT,HZ,TZK} and the references therein). This kind of BSDEs,
which are reflected along an oblique reflection rather than a normal one in a convex domain, were first introduced by  Hamadène and Jeanblanc \cite{HJ}, where they used its solution to characterize the value of an optimal switching problem, in particular in the setting of power plant management. The related equation was solved by Hu and Tang \cite{HT} using the penalization method and by Hamadène and Zhang \cite{HZ} using the Picard iteration method, they generalized the preceding work. See also Chassagneux et al. \cite{CEK}. Recently, Chassagneux and Richou \cite{CR} obtained the existence and uniqueness to multidimensional RBSDEs in an open convex domain, allowing for oblique directions of reflection. And then, B\'en\'ezet et al. \cite{BCR} introduced and studied a new class of optimal switching problems, namely controlled randomization problems, where some extra-randomness impacts the choice of switching modes and associated costs.
\par In the case when the horizon is infinite, there is still much to be done and this
is the novelty of this paper. So, our problem is how to generalize the multi-dimensional RBSDEs to an infinite horizon. In that case, El Asri \cite{Elasri} studied RBSDEs for which the generator does not depend on $(y,z)$ and provided an application to optimal switching problem. The latter is a problem in which a decision controller controls a system which may operate in different modes (e.g., a power plant). The aim of the controller is to maximize some performance criterion by optimally choosing controls of the form $a:=(\tau_{n},\zeta_{n})_{n\geq0}$. Here $(\tau_{n})_{n\geq0}$ denotes an increasing sequence of (stopping) times at which the controller switches the system across different operating modes. Moreover, $(\zeta_{n})_{n\geq0}$ is a sequence of random variables taking their values in $\I$. Each $\zeta_{n}$ represents the system's new operating mode after a switch has occurred at time $\tau_{n}$. In this setting it is well known that $Y^{i}_t$ is the value of an optimal switching strategy, i.e.,
\begin{equation}
e^{-rt}Y^{i}_t = \underset{(\tau_{n},\zeta_{n})_{n\geq0}\in\D_{t}^{i}}{\esssup}\;\E\biggl[\int_{t}^{+\infty} e^{-rs}f_{a_s}(s,X_s) \, ds-A_{\infty}^{a}\biggr];
\end{equation}
where $\D_{t}^{i}=\{(\tau_{n},\zeta_{n})_{n\geq0}  \mbox{ such that } \tau_{0}=t \mbox{ and } \zeta_{0}=i\}$, the process $(a_s)_{s\geq0}$ is indicating the mode of the system at time $s$ and $A_{\infty}^{a}$ stands for the total switching cost when the strategy $a$ is implemented.
\par The main contribution of this paper is to establish the existence and uniqueness of solution for RBSDE $(\ref{rbsde})$. First, we prove the existence and uniqueness (when, the generator $f_i$ does not depend on all $Y^i$) by penalization and verification methods, respectively. And then, using the contraction method, we obtain the existence and uniqueness result for RBSDE $(\ref{rbsde})$.\\
\par The paper is organized as follows: In section $2$, we state some notations and assumptions. In section $3$, we prove the existence and uniqueness of solution when the generator does depend only on $Y^i$. Finally, in Section $4
$, we state and prove the main result concerning the existence and uniqueness of  solutions to the system of  RBSDE $(\ref{rbsde})$.
\section{Notations and Assumptions}
Throughout this paper, let $(\Omega,\mathcal{F},\P)$ is a fixed probability space endowed with a $d$-dimensional Brownian motion $B=(B_t)_{t\geq0}$, where $\{\mathcal{F}_t,\, t\geq0\}$ is the natural filtration of the Brownian motion augmented by $\P$-null
sets of $\mathcal{F}$, with $\F_{\infty}=\bigvee_{t\geq0}\F_t$. All the measurability notions will refer to this filtration. Let $|.|$ denote the Euclidean norm for vectors.\\
\par Let us introduce the following spaces:
\begin{itemize}
\item[-] $\p$ is the $\sigma$-algebra on $[0,+\infty)\times\Omega$ of $\F$-progressively measurable sets.
	\item[-]$\mathcal{M}^{2}$  is the set of $\R^{d}$-valued, progressively measurable processes $(Z_t)_{t\geq0}$ such that $$\mathbb{E}\left[\int^{+\infty}_{0}|Z_s|^2ds\right]<+\infty.$$
	\item[-]$\mathcal{S}^{2}$ is the set of $\mathbb{R}$-valued adapted and  c\`adl\`ag processes $(Y_t)_{t\geq0}$ such that $$\mathbb{E}\left[\sup_{t\geq0}|Y_{t}|^2\right]<+\infty.$$
	\item[-]$\mathcal{K}^2$  is the subset of non-decreasing processes $(K_t)_{t\geq 0}\in\mathcal{S}^{2}$, starting from
	 $K_0=0$.	
	\item[-] $\L^2$  is the
	set of $\F_{\infty}$-measurable random 	variables $\xi$ satisfying $\E[|\xi|^2]<+\infty$.
\end{itemize}
Let us now consider the following function: for $i,j\in\mathcal{I}$,
	$$f_i(t,x,y^1,...,y^m,z^i):\Omega\times[0,+\infty)\times\R^k\times\R^{m}\times\R^d\mapsto\R.$$
Next, let us consider a deterministic and bounded function $u(t):[0,+\infty)\mapsto[0,+\infty)$ which satisfies: $\int_0^{+\infty}u(t)dt<+\infty$ and $\int_0^{+\infty}u^2(t)dt<+\infty$.  In what follows we take for simplicity the same function $u(t)$ in the inequalities below, in addition we denote $\vec{y}=(y^1,...,y^m).$ \\
Now, we make the following assumptions:
\begin{itemize}
\item[$\mathbf{[H1]}$] For $i\in\mathcal{I}$, $f_i:\Omega\times[0,+\infty)\times\R^k\times\R^{m}\times\R^d\mapsto\R$ is $\mathcal{B}([0,+\infty))\otimes \mathcal{B}(\R^k)$-measurable and satisfies:\\
$(i)$ for any $i\in\I$, $t\in[0,+\infty)$ and $x\in\R^k$, $f_i(t,x,0,0)$ belongs to $\mathcal{M}^2$.\\
$(ii)$ $f_i(t,x,\vec{y},z)$ is Lipschitz continuous in $(\vec{y},z),$ i.e., for all $(t,x,\vec{y}_j,z_j)\in[0,+\infty)\times\R^{k}\times \R^{m}\times \R^{d},$ $j=1,2
$  we have
\begin{equation}
|f_i(t,x,\vec{y}_1,z_1)-f_i(t,x,\vec{y}_2,z_2)|\leq u(t)(|\vec{y}_1-\vec{y}_2|+|z_1-z_2|).
\end{equation}
\item[$\mathbf{[H2]}$] For any $(i,j)\in\I^2$, $g_{ij}$ satisfies the following:
\begin{itemize}
	\item [(i)] $g_{ii}=0;$
	\item [(ii)] $g_{ij}>0,$ for $i\neq j$;
	\item [(iii)] for any $(i,j,l)\in\mathcal{I}^3$, such that $i\neq j$ and $j\neq l$, we have $$g_{ij}+g_{jl}\geq g_{il}.$$
\end{itemize}
\item[$\mathbf{[H3]}$] For any $i\in\mathcal{I}$, $\xi^i:=\underset{t \to+\infty}{\lim}Y^i_t$ belongs to $\L^2$ and satisfies $$\xi^i\geq \underset{j\in\mathcal{I}^{-i}}{\max}(\xi^j-g_{ij}).$$
\end{itemize}
\begin{remark}
	Notice that in the event that $g_{ij}$ is not constant and depends on $t$ and the process $X$. To show the main result of this paper, we shall require the following assumption:\\
	$\mathbf{[Hg]}$ For any $(i,j)\in\I^2$, $g_{ij}$ belongs to $\mathcal{C}^{1,2}([0,+\infty)\times\R^{k})$ and $$(\partial_tg_{ij}+\mathcal{L}_x g_{ij})(t,x)\leq0, \mbox{ for all } (t,x)\in[0,+\infty)\times\R^{k},$$ where $\mathcal{L}_x$ is the infinitesimal generator of $X$.
\end{remark}
Before solving RBSDE $(\ref{rbsde})$, we first provide its existence and uniqueness when, for any $i\in\I$, the function $f_i$ does not depend on all $Y^j,\;\forall j\in\I^{-i}$, that is, $\P$-a.s., $f_i(t,x,\vec{y},z^i):= f_i(t,x,y^i,z^i),$ for any $t,x,y^i,\;and\;z^i$, and consider the following RBSDE, $\forall t\geq0$
\begin{equation}\label{keybsde}
\begin{cases}
e^{-rt}Y^{i}_t = \int_{t}^{+\infty} e^{-rs}f_i(s,X_s,Y^i_s,Z^{i}_s) \, ds + K^{i}_{\infty} - K^{i}_t-  \int_{t}^{+\infty}e^{-rs}Z^{i}_s \, dB_s,\\
\forall t \geq 0, \quad e^{-rt}Y^{i}_t \geq e^{-rt} \underset{j\in\mathcal{I}^{-i}}{\max}(Y^{j}_t-g_{ij}),\\
\int_{0}^{+\infty} e^{-rs}\{Y^{i}_s - \underset{j\in\mathcal{I}^{-i}}{\max}(Y^{j}_s-g_{ij})\} \, dK^i_s = 0.
\end{cases}
\end{equation}
\section{Existence and Uniqueness of a Solution for RBSDE $(\ref{keybsde})$}
\subsection{Existence}
In this subsection, we shall prove an existence result of the solution of RBSDE $(\ref{keybsde}).$
\begin{proposition}
	Under $\mathbf{[H1]}$, $\mathbf{[H2]}$ and $\mathbf{[H3]}$, the RBSDE $(\ref{keybsde})$ has at least one solution.
\end{proposition}\label{prop1}
\begin{proof}
	The proof will be divided into four steps. \\
	\\ \underline{\bf{Step 1}}: The penalized BSDE.\\
	For any $i\in\mathcal{I}$ and $n\geq1$  let us consider $(Y^{i,n},Z^{i,n})\in\mathcal{S}^2\times\mathcal{M}^2$ the unique solution of the following BSDE:
	\begin{equation}\label{pbsde}
	e^{-rt}Y^{i,n}_t = \int_{t}^{+\infty} e^{-rs}f_i^n(s,X_s,Y^{i,n}_s,Z^{i,n}_s) \, ds-  \int_{t}^{+\infty}e^{-rs}Z^{i,n}_s \, dB_s;
	\end{equation}
	where $f_i^n$ is defined on $[0,+\infty)\times\R^k\times\R\times\R^d$ by $$f_i^n:(t,x,y^i,z^i)\mapsto f_i(t,x,y^i,z^i) +n\sum_{j=1}^{m}(y^{i}-y^{j}+g_{ij})^-.$$
	Indeed, thanks to the result by Chen \cite{chen}, this solution exists and is unique.\\
	\\ \underline{\bf{Step 2}}: A priori estimate.\\
	In this step, we derive two lemmas on the a priori estimation of the penalized BSDE $(\ref{pbsde})$, which will play a primordial role in the proof of Proposition $(\ref{prop1})$.
	\begin{lemma}\label{lemma31}
		There exist a positive constant $C_u$ which depends on $u$, such that the following hold true: for any $i,j\in\mathcal{I}$ and $n\geq1$.
		\begin{equation}\label{Yijnestime}
		\begin{aligned}
		\E&\biggl[\underset{t\geq0}{\sup}\,[(Y_t^{i,n}-Y_t^{j,n}+g_{ij})^-]^2+n^2\int_{0}^{+\infty}[(Y_s^{i,n}-Y_s^{j,n}+g_{ij})^-]^2ds\biggr]\\ &\leq \E\biggl[\int_{0}^{+\infty}C_u1_{\L_{ij,n}}(s)[1+|f_i(s,X_s,0,0)|^2+|f_j(s,X_s,0,0)|^2+|Y_s^{i,n}|^2+|Z_s^{i,n}|^2]ds\biggr],
		\end{aligned}
		\end{equation}
		where $$\L_{ij,n}:=\{(s,\omega)\in[0,+\infty)\times\Omega, \mbox{ such that } Y_s^{i,n}-Y_s^{j,n}+g_{ij}<0 \}.$$
	\end{lemma}
	\begin{proof} For given $i,j\in\I$, set
		\begin{equation}\label{yijn}
			\Y^{ij,n}_t:=Y_t^{i,n}-Y_t^{j,n}+g_{ij}, \qquad t\in[0,+\infty).
		\end{equation}
		Next, rewrite equation  $(\ref{pbsde})$ in differential form
		\begin{equation}\label{dpbsde}
		d(e^{-rt}Y^{i,n}_t) = -e^{-rt}f_i^n(t,X_t,Y_t^{i,n},Z^{i,n}_t) \, dt +e^{-rt}Z^{i,n}_t \, dB_t.
		\end{equation}
		So for $t\geq0$, $n\geq1$ and $i,l\in\I$, equation $(\ref{pbsde})$ is equivalent to
		\begin{equation}\label{yin}
		Y^{i,n}_t = \xi^i+\int_{t}^{+\infty} [f_i(s,X_s,Y^{i,n}_s,Z^{i,n}_s) +n\sum_{l=1}^{m}(\Y^{il,n}_s)^- -rY^{i,n}_s]\, ds-  \int_{t}^{+\infty}Z^{i,n}_s \, dB_s;
		\end{equation}
		where $\xi^i:=\underset{t \to+\infty}{\lim}Y^{i,n}_{t},\; \forall i\in\I$ (this limit exists thanks to Assumption $\mathbf{[H3]}$).\\
		By an application of Itô-Tanaka's formula, for every
		$t\geq0$, we obtain
		\begin{equation}\label{ito-tanaka}
		\begin{aligned}
		d([(\Y^{ij,n}_t)^-]^2)=&-2(\Y^{ij,n}_t)^-d\Y^{ij,n}_t+\frac{1}{2}\Y^{ij,n}_tdL_t^0(\Y^{ij,n})\\ &+1_{\L_{ij,n}}(t)(Z^{i,n}_t-Z^{j,n}_t )^2dt;
		\end{aligned}
		\end{equation}
		 where $L^0(\Y^{ij,n})$ denotes the local-time at zero of the semi-martingale $\Y^{ij,n}$.\\
		We now want to find a convenient expression for $d\Y^{ij,n}_t$. In the definition of $\Y^{ij,n}_t$ (cf.
		$(\ref{yijn})$) we may express $Y^{i,n}$ and $Y^{j,n}$
		in terms of their associated BSDE $(\ref{yin})$.  This
		gives, for any $t\in[0,+\infty)$
		\begin{equation}\label{Yijn}
		\begin{aligned}
		\Y^{ij,n}_t =& \xi^i-\xi^j+g_{ij}+\int_{t}^{+\infty} [f_i(s,X_s,Y^{i,n}_s,Z^{i,n}_s)-f_j(s,X_s,Y^{j,n}_s,Z^{j,n}_s)\\ &-r(Y^{i,n}_s-Y^{j,n}_s)]\, ds +n\int_{t}^{+\infty}\sum_{l=1}^{m}(\Y^{il,n}_s)^-\, ds-n\int_{t}^{+\infty}\sum_{l=1}^{m}(\Y^{jl,n}_s)^-\, ds\\ &-\int_{t}^{+\infty}(Z^{i,n}_s-Z^{j,n}_s) \, dB_s.
		\end{aligned}
		\end{equation}
		Then,
		\begin{equation}\label{dyijn}
		\begin{aligned}
		d\Y^{ij,n}_t =& -[f_i(t,X_t,Y^{i,n}_t,Z^{i,n}_t)-f_j(t,X_t,Y^{j,n}_t,Z^{j,n}_t)-r(Y^{i,n}_t-Y^{j,n}_t)]\, dt\\ & -n\sum_{l=1}^{m}(\Y^{il,n}_t)^-\, dt+n\sum_{l=1}^{m}(\Y^{jl,n}_t)^-\, dt+(Z^{i,n}_t-Z^{j,n}_t )\, dB_t.
    	\end{aligned}
\end{equation}
Noticing that the integral with respect to the local-time $L^0(\Y^{ij,n})$ is zero and by Assumption $\mathbf{[H3]}$ we have that $(\xi^i-\xi^j+g_{ij})^-=0$,  we obtain from $(\ref{ito-tanaka})$ that for
every $t\in[0,+\infty)$
\begin{equation}\label{tanaka}
\begin{aligned}
&[(\Y_t^{ij,n})^-]^2+2n\int_{t}^{+\infty}[(\Y_s^{ij,n})^-]^2ds +\int_{t}^{+\infty}1_{\L_{ij,n}}(s)(Z_s^{i,n}-Z_s^{j,n})^2ds\\ &= 2\int_{t}^{+\infty}(\Y_s^{ij,n})^-[f_j(s,X_s,Y_s^{j,n},Z_s^{j,n})-f_i(s,X_s,Y_s^{i,n},Z_s^{i,n})+r(Y_s^{i,n}-Y_s^{j,n})]ds \\ & \quad+2\int_{t}^{+\infty}(\Y_s^{ij,n})^-(Z_s^{i,n}-Z_s^{j,n})dB_s +2n\int_{t}^{+\infty}(\Y_s^{ij,n})^-(\Y_s^{ji,n})^-ds\\ &\quad+2n\underset{l\neq i,j}{\sum}\int_{t}^{+\infty}(\Y_s^{ij,n})^-[(\Y_s^{jl,n})^--(\Y_s^{il,n})^-]ds.
\end{aligned}
\end{equation}	
But, by Assumption $\mathbf{[H2]}-(iii)$ we have that $g_{ij}+g_{ji}\geq g_{ii}= 0$. Thus, we obtain that, for every $s\in[0,+\infty)$
$$\{y\in\R^m,\,\, y^i_s-y^j_s+g_{ij}<0\}\cap\{y\in\R^m,\,\,y^j_s-y^i_s+g_{ji}<0\}=\emptyset .$$
from which we deduce that $$(\Y_s^{ij,n})^-(\Y_s^{ji,n})^-=0,\qquad i,j\in\I.$$
Relying next on the elementary inequality $x_1^--x_2^- \leq (x_1-x_2)^-$ and by Assumption $\mathbf{[H2]}-(iii)$, we get
\begin{equation*}
\begin{aligned}
 (\Y_s^{ij,n})^-[(\Y_s^{jl,n})^--(\Y_s^{il,n})^-] &\leq (Y_s^{i,n}-Y_s^{j,n}+g_{ij})^-(Y_s^{j,n}-Y_s^{i,n}+g_{jl}-g_{il})^-\\ & = (Y_s^{i,n}-Y_s^{j,n}+g_{ij})^-(Y_s^{i,n}-Y_s^{j,n}+g_{il}-g_{jl})^+ \\ & \leq (Y_s^{i,n}-Y_s^{j,n}+g_{ij})^-(Y_s^{i,n}-Y_s^{j,n}+g_{ij})^+ \\ & =0.
\end{aligned}
\end{equation*}
Then, taking expectation on both sides of $(\ref{tanaka}),$ we obtain
\begin{equation}\label{Exp}
\begin{aligned}
&\E\biggl[[(\Y_t^{ij,n})^-]^2\biggr]+2n\E\biggl[\int_{t}^{+\infty}[(\Y_s^{ij,n})^-]^2ds\biggr]+\E\biggl[\int_{t}^{+\infty}1_{\mathcal{L}_{ij,n}}(s)|Z_s^{i,n}-Z_s^{j,n}|^2ds\biggr]\\ &\leq  2\E\biggl[\int_{t}^{+\infty}(\Y_s^{ij,n})^-\{|f_j(s,X_s,Y_s^{j,n},Z_s^{j,n})-f_i(s,X_s,Y_s^{i,n},Z_s^{i,n})|+r|Y_s^{i,n}-Y_s^{j,n}|\}ds\biggr].
\end{aligned}
\end{equation}
	Noting that,
\begin{equation}\label{fi-fj}
\begin{aligned}
&	|f_j(s,X_s,Y_s^{j,n},Z_s^{j,n})-f_i(s,X_s,Y_s^{i,n},Z_s^{i,n})| +r|Y_s^{i,n}-Y_s^{j,n}|\\ & \leq |f_i(s,X_s,Y_s^{i,n},Z_s^{i,n})-f_j(s,X_s,Y_s^{i,n},Z_s^{i,n})|+|f_j(s,X_s,Y_s^{i,n},Z_s^{i,n})-f_j(s,X_s,Y_s^{j,n},Z_s^{j,n})|\\ & \qquad +r|Y_s^{i,n}-Y_s^{j,n}|,\\ & \leq |f_i(s,X_s,Y_s^{i,n},Z_s^{i,n})|+|f_j(s,X_s,Y_s^{i,n},Z_s^{i,n})|+|f_j(s,X_s,Y_s^{i,n},Z_s^{i,n})-f_j(s,X_s,Y_s^{j,n},Z_s^{j,n})|\\ & \qquad+r|Y_s^{i,n}-Y_s^{j,n}|, \\ & \leq |f_i(s,X_s,0,0)|+|f_j(s,X_s,0,0)|+(r+u(s))|Y_s^{i,n}-Y_s^{j,n}|+u(s)[2|Y_s^{i,n}|+2|Z_s^{i,n}|\\ &\qquad +|Z_s^{i,n}-Z_s^{j,n}|] ,\\ & \leq C_u[1+|f_i(s,X_s,0,0)|+|f_j(s,X_s,0,0)|+|Y_s^{i,n}|+|Z_s^{i,n}|+|\Y_s^{ij,n}|+|Z_s^{i,n}-Z_s^{j,n}|];
\end{aligned}
\end{equation}
where $C_u$ is a constant depending on $u(s)$, independent of $n$ and which might hereafter vary from line to line.\\ Now going back to $(\ref{Exp})$ and using $(\ref{fi-fj}),$ we get that
\begin{equation}\label{Expp}
\begin{aligned}
&\E\biggl[[(\Y_t^{ij,n})^-]^2\biggr]+2n\E\biggl[\int_{t}^{+\infty}[(\Y_s^{ij,n})^-]^2ds\biggr]+\E\biggl[\int_{t}^{+\infty}1_{\L_{ij,n}}(s)|Z_s^{i,n}-Z_s^{j,n}|^2ds\biggr]\\ &\leq  2\E\biggl[\int_{t}^{+\infty}C_u(\Y_s^{ij,n})^-[1+|f_i(s,X_s,0,0)|+|f_j(s,X_s,0,0)|+|Y_s^{i,n}|+|Z_s^{i,n}|+|\Y_s^{ij,n}|\\ & \qquad\quad+|Z_s^{i,n}-Z_s^{j,n}|]ds\biggr],\\ &\leq  \E\biggl[\int_{t}^{+\infty}C_u[(\Y_s^{ij,n})^-]^2ds\biggr]+
\frac{1}{2}\E\biggl[\int_{t}^{+\infty}1_{\L_{ij,n}}(s)[1+|f_i(s,X_s,0,0)|^2+|f_j(s,X_s,0,0)|^2\\ & \qquad\quad+|Y_s^{i,n}|^2+|Z_s^{i,n}|^2+|(\Y_s^{ij,n})^-|^2+|Z_s^{i,n}-Z_s^{j,n}|^2]ds\biggr].
\end{aligned}
\end{equation}
Applying Gronwall's inequality, it follows that
\begin{equation}\label{gronwal}
\E\biggl[[(\Y_t^{ij,n})^-]^2\biggr]\leq \E\biggl[\int_{0}^{+\infty}C_u1_{\L_{ij,n}}(s)[1+|f_i(s,X_s,0,0)|^2+|f_j(s,X_s,0,0)|^2+|Y_s^{i,n}|^2+|Z_s^{i,n}|^2]ds\biggr],
\end{equation}
and
\begin{equation}\label{nEyij-}
\begin{aligned}
 n\E\biggl[&\int_{0}^{+\infty}[(\Y_s^{ij,n})^-]^2ds\biggr]+\E\biggl[\int_{0}^{+\infty}1_{\L_{ij,n}}(s)|Z_s^{i,n}-Z_s^{j,n}|^2ds\biggr]\\ &\leq
\E\biggl[\int_{0}^{+\infty}C_u1_{\L_{ij,n}}(s)[1+|f_i(s,X_s,0,0)|^2+|f_j(s,X_s,0,0)|^2+|Y_s^{i,n}|^2+|Z_s^{i,n}|^2]ds\biggr].
\end{aligned}
\end{equation}
Going back to $(\ref{tanaka})$ and applying Burkholder-Davis-Gundy's inequality, we obtain
\begin{equation}\label{supBDG}
\begin{aligned}
\E\biggl[\underset{t\geq0}{\sup}\,[(\Y_t^{ij,n})^-]^2\biggr]\leq& \E\biggl[\int_{0}^{+\infty}C_u1_{\L_{ij,n}}(s)[1+|f_i(s,X_s,0,0)|^2+|f_j(s,X_s,0,0)|^2\\ & \qquad\quad+|Y_s^{i,n}|^2+|Z_s^{i,n}|^2]ds\biggr].
\end{aligned}
\end{equation}
On the other hand, from $(\ref{Exp})$, we deduce that,
\begin{equation*}
\begin{aligned}
2n&\E\biggl[\int_{0}^{+\infty}[(\Y_s^{ij,n})^-]^2ds\biggr]\\ &\leq  \E\biggl[\int_{0}^{+\infty}(n+C_u)[(\Y_s^{ij,n})^-]^2ds\biggr]+
\frac{1}{n}\E\biggl[\int_{0}^{+\infty}C_u1_{\L_{ij,n}}(s)[1+|f_i(s,X_s,0,0)|^2+|f_j(s,X_s,0,0)|^2\\ & \qquad\quad+|Y_s^{i,n}|^2+|Z_s^{i,n}|^2]ds\biggr].
\end{aligned}
\end{equation*}
For $n$ large enough we finally deduce that
\begin{equation}\label{fin}
\begin{aligned}
n^2\E\biggl[\int_{0}^{+\infty}[(\Y_s^{ij,n})^-]^2ds\biggr]\leq& \E\biggl[\int_{0}^{+\infty}C_u1_{\L_{ij,n}}(s)[1+|f_i(s,X_s,0,0)|^2+|f_j(s,X_s,0,0)|^2\\ & \qquad\quad+|Y_s^{i,n}|^2+|Z_s^{i,n}|^2]ds\biggr].
\end{aligned}
\end{equation}
\end{proof}
Thanks to Lemma $(\ref{lemma31})$ we are able to prove the next uniform estimate on the solution of the penalized problem.
\begin{lemma}\label{lemma32}  For any $i,j\in\mathcal{I},$ $n\geq1.$ There exist a positive constant $C$ independent of $n$ such that,
	\begin{equation}\label{Ynestime}
	 \E\biggl[\underset{t\geq0}{\sup}\,|Y_t^{i,n}|^2+\int_{0}^{+\infty}|Z_s^{i,n}|^2ds+n^2\int_{0}^{+\infty}\sum_{j=1}^{m}[(Y_s^{i,n}-Y_s^{j,n}+g_{ij})^-]^2ds\biggr]\leq C.
	\end{equation}
\end{lemma}
\begin{proof}
	Applying Itô's formula with $|Y_t^{i,n}|^2$ and recalling $(\ref{yin})$ we obtain that for any $t\in[0,+\infty)$
		\begin{equation}\label{ito_exp_y}
	\begin{aligned}
	|Y_t^{i,n}|^2+&\int_{t}^{+\infty}[2r|Y_s^{i,n}|^2+|Z_s^{i,n}|^2]ds = 2\int_{t}^{+\infty}Y_s^{i,n}[f_i(s,X_s,Y^{i,n}_s,Z_s^{i,n}) \\ & +n\sum_{j=1}^{m}(Y^{i,n}_s-Y^{j,n}_s+g_{ij})^-]ds -2\int_{t}^{+\infty}Y_s^{i,n}Z_s^{i,n}dB_s.
	\end{aligned}
	\end{equation}
	Taking expectation and using  the classical inequality $2ab\leq \frac{1}{\epsilon}a^2+\epsilon b^2,$ for any $a,b\in\R$ and $\epsilon>0$, we obtain
	\begin{equation*}
	\begin{aligned}
	& \E[|Y_t^{i,n}|^2]+\E\biggl[\int_{t}^{+\infty} (2r|Y_s^{i,n}|^2+|Z_s^{i,n}|^2)ds\biggr] \\ & \leq 2\E\biggl[\int_{t}^{+\infty}|Y_s^{i,n}|[|f_i(s,X_s,Y^{i,n}_s,Z_s^{i,n})|  +n\sum_{j=1}^{m}(Y^{i,n}_s-Y^{j,n}_s+g_{ij})^-]ds\biggr], \\ & \leq 2\E\biggl[\int_{t}^{+\infty}|Y_s^{i,n}|[|f_i(s,X_s,0,0)|+u(s)(|Y_s^{i,n}|+|Z_s^{i,n}| ) \\ & \qquad +n\sum_{j=1}^{m}(Y^{i,n}_s-Y^{j,n}_s+g_{ij})^-]ds\biggr], \\ & \leq \E\biggl[\int_{t}^{+\infty}[2(u^2(s)+u(s)+\epsilon)|Y_s^{i,n}|^2+\frac{1}{2}|Z_s^{i,n}|^2+\epsilon^{-1}|f_i(s,X_s,0,0)|^2]ds\biggr] \\ &  \qquad + \epsilon^{-1} n^2 \E\biggl[\int_{t}^{+\infty}\sum_{j=1}^{m}[(Y^{i,n}_s-Y^{j,n}_s+g_{ij})^-]^2ds\biggr].
	\end{aligned}
	\end{equation*}
	In view of $(\ref{Yijnestime})$, if we choose $\epsilon=4C_u$ and for $r \geq  u^2(s)+u(s)+4C_u+\frac{1}{8}$, we obtain
	\begin{equation}
	\E[|Y_t^{i,n}|^2]+\frac{1}{4}\E\biggl[\int_{t}^{+\infty}|Z_s^{i,n}|^2ds\biggr]\leq \E\biggl[\int_{t}^{+\infty}C_u[1+|f_i(s,X_s,0,0)|^2+|f_j(s,X_s,0,0)|^2]ds\biggr].
	\end{equation}
	Therfore, for $t=0$ we have
	\begin{equation}
	\E\big[\int_{0}^{+\infty}|Z_s^{i,n}|^2ds\big]\leq 4\E\biggl[\int_{0}^{+\infty}C_u[1+|f_i(s,X_s,0,0)|^2+|f_j(s,X_s,0,0)|^2]ds\biggr].
	\end{equation}
	Using again equation $(\ref{ito_exp_y})$ and the Burkholder-Davis-Gundy's inequality, for some finite universal constant $c$, we obtain that
	\begin{equation}
	\begin{aligned}
	\E[\underset{t\geq0}{\sup}\,|Y_t^{i,n}|^2]\leq & \; \E\biggl[\int_{t}^{+\infty}C_u[1+|f_i(s,X_s,0,0)|^2+|f_j(s,X_s,0,0)|^2]ds\biggr]\\ &\;+c\E\biggl[(\int_{0}^{+\infty}|Y_s^{i,n}|^2|Z_s^{i,n}|^2ds)^{\frac{1}{2}}\biggr];
	\end{aligned}
	\end{equation}
	with,
	\begin{equation}
	\begin{aligned}
	c\E\biggl[\big(\int_{0}^{+\infty}|Y_s^{i,n}|^2|Z_s^{i,n}|^2ds\big)^{\frac{1}{2}}\biggr]\leq\; \frac{1}{2}\E\biggl[\underset{t\geq0}{\sup}\,|Y_t^{i,n}|^2\biggr]+\frac{c^2}{2}\E\biggl[\int_{0}^{+\infty}|Z_s^{i,n}|^2ds\biggr].
	\end{aligned}
	\end{equation}
	Then,
	\begin{equation}
	\begin{aligned}
	\E[\underset{t\geq0}{\sup}\,|Y_t^{i,n}|^2]\leq \; 2\E\biggl[\int_{t}^{+\infty}&C_u[1+|f_i(s,X_s,0,0)|^2+|f_j(s,X_s,0,0)|^2]ds\biggr]\\ &\;+c^2\E[\int_{0}^{+\infty}|Z_s^{i,n}|^2ds].
	\end{aligned}
	\end{equation}
	Therefore,
	\begin{equation}
	\E[\underset{t\geq0}{\sup}\,|Y_t^{i,n}|^2]+\E[\int_{0}^{+\infty}|Z_s^{i,n}|^2ds]\leq C.
	\end{equation}
	From $(\ref{Yijnestime})$, we obtain
	\begin{equation*}
	n^2\E\biggl[\int_{0}^{+\infty}[(Y_s^{i,n}-Y_s^{j,n}+g_{ij})^-]^2ds\biggr]\leq C.
	\end{equation*}
	Taking the summation over all $j\in\I$, we obtain
	\begin{equation}\label{n2}
	n^2\E\biggl[\int_{0}^{+\infty}\sum_{j=1}^{m}[(Y_s^{i,n}-Y_s^{j,n}+g_{ij})^-]^2ds\biggr]\leq C.
	\end{equation}
	The proof of Lemma $(\ref{lemma32})$ is now complete.
\end{proof}
\underline{\bf{Step 3}}: Convergence of the sequence.\\
In order to show that the sequence $(Y_t^{i,n})_{n\geq1}$ is non-decreasing and convergent for any $i\in\mathcal{I}$ we use a comparison theorem for infinite horizon BSDE presented in \cite{HLW}.
Since $f_i^n\leq f_i^{n+1}$, by the comparison Theorem, we have, for any $n\geq1$  $$Y_t^{i,n}\leq Y_t^{i,n+1},\mbox{ for all } i\in\mathcal{I} \mbox{ and } t\in[0,+\infty).$$
Then, $Y_t^{i,n}$ admits a limit denoted by $Y_t^{i}$. Moreover, from the a priori  estimate $(\ref{Ynestime})$ and Fatou's Lemma, we have
\begin{equation}\label{Yestime}
\E[\underset{t\geq0}{\sup}\,|Y_t^{i}|^2]\leq C.
\end{equation}
Then, applying Lebesgue's dominated converge theorem, we get
\begin{equation}\label{lebesgue}
\underset{n\to+\infty}{\lim}\E\biggl[\int_{0}^{+\infty}|Y_t^{i,n}-Y_t^{i}|^2ds\biggr]=0.
\end{equation}
Now, we prove that $(Y^{i,n},Z^{i,n})$ is a Cauchy sequence. To do so we apply Itô's formula to $|Y_t^{i,n}-Y_t^{i,p}|^2,$ to get
\begin{equation}\label{ynyp2}
\begin{aligned}
&|Y_t^{i,n}-Y_t^{i,p}|^2+\int_{t}^{+\infty}(2r |Y_s^{i,n}-Y_s^{i,p}|^2+|Z_s^{i,n}-Z_s^{i,p}|^2)ds\\ &= 2\int_{t}^{+\infty}(Y_s^{i,n}-Y_s^{i,p})(f_i(s,X_s,Y_s^{i,n},Z_s^{i,n})-f_i(s,X_s,Y^{i,p}_s,Z_s^{i,p}))ds\\ & \qquad+2n\sum_{j=1}^{m}\int_{t}^{+\infty}(Y_s^{i,n}-Y_s^{i,p})(Y^{i,n}_s-Y^{j,n}_s+g_{ij})^-ds\\&\qquad-2p\sum_{j=1}^{m}\int_{t}^{+\infty}(Y_s^{i,n}-Y_s^{i,p})(Y^{i,p}_s-Y^{j,p}_s+g_{ij})^-ds\\&\qquad-2\int_{t}^{+\infty}(Y_s^{i,n}-Y_s^{i,p})(Z_s^{i,n}-Z_s^{i,p})dB_s.
\end{aligned}
\end{equation}
Taking expectation on both sides of the last equality yields
\begin{equation*}
\begin{aligned}
&\E[|Y_t^{i,n}-Y_t^{i,p}|^2]+\E\biggl[\int_{t}^{+\infty}(2r |Y_s^{i,n}-Y_s^{i,p}|^2+|Z_s^{i,n}-Z_s^{i,p}|^2)ds\biggr]\\ &\leq 2\E\biggl[\int_{t}^{+\infty}|Y_s^{i,n}-Y_s^{i,p}||f_i(s,X_s,Y_s^{i,n},Z_s^{i,n})-f_i(s,X_s,Y^{i,p}_s,Z_s^{i,p})|ds\biggr]\\ & \qquad+2n\sum_{j=1}^{m}\E\biggl[\int_{t}^{+\infty}|Y_s^{i,n}-Y_s^{i,p}|(Y^{i,n}_s-Y^{j,n}_s+g_{ij})^-ds\biggr]\\ & \qquad+2p\sum_{j=1}^{m}\E\biggl[\int_{t}^{+\infty}|Y_s^{i,n}-Y_s^{i,p}|(Y^{i,p}_s-Y^{j,p}_s+g_{ij})^-ds\biggr], \\ & \leq \E\biggl[\int_{t}^{+\infty}2(u^2(s)+u(s))|Y_s^{i,n}-Y_s^{i,p}|^2ds\biggr]+\frac{1}{2}\E\biggl[\int_{t}^{+\infty}|Z_s^{i,n}-Z_s^{i,p}|^2ds\biggr]\\ & \qquad +2\E\biggl[\int_{t}^{+\infty}|Y_s^{i,n}-Y_s^{i,p}|^2ds\biggr]^{\frac{1}{2}}\sum_{j=1}^{m}\biggl(n^2\E\biggl[\int_{t}^{+\infty}[(Y^{i,n}_s-Y^{j,n}_s+g_{ij})^{-}]^2ds\biggr]\biggr)^{\frac{1}{2}}\\ &\qquad +2\E\biggl[\int_{t}^{+\infty}|Y_s^{i,n}-Y_s^{i,p}|^2ds\biggr]^{\frac{1}{2}}\sum_{j=1}^{m}\biggl(p^2\E\biggl[\int_{t}^{+\infty}[(Y^{i,p}_s-Y^{j,p}_s+g_{ij})^{-}]^2ds\biggr]\biggr)^{\frac{1}{2}}.
\end{aligned}
\end{equation*}
Setting $t=0,$ and choosing $r \geq u^2(s)+u(s)+1$, from $(\ref{Ynestime})$ and $(\ref{lebesgue}) $, we have that
\begin{equation}\label{Znp}
\underset{n,\, p \to +\infty}{lim}\mathbb{E}\left[\int_{0}^{+\infty}|Z_s^{i,n}-Z_s^{i,p}|^2 \, \mathrm{d}{s}\right]=0.
\end{equation}
Consequently, the sequence $(Z^{i,n})_{n\geq1}$ converges in $\mathcal{M}^2$ to a process which we denote $Z^i$.\\
Now, going back to $(\ref{ynyp2}) $, we deduce that
\begin{equation}\label{BDG}
\begin{aligned}
\E\biggl[\underset{t\geq0}{\sup}\,|Y_t^{i,n}-Y_t^{i,p}|^2\biggr]&\leq 2\E\biggl[\int_{t}^{+\infty}|Y_s^{i,n}-Y_s^{i,p}||f_i(s,X_s,Y_s^{i,n},Z_s^{i,n})-f_i(s,X_s,Y^{i,p}_s,Z_s^{i,p})|ds\biggr]\\ &\qquad +2n\sum_{l=1}^{m}\E\biggl[\int_{t}^{+\infty}|Y_s^{i,n}-Y_s^{i,p}|(Y^{i,n}_s-Y^{l,n}_s+g_{il})^-ds\biggr]\\ &\qquad +2p\sum_{l=1}^{m}\E\biggl[\int_{t}^{+\infty}|Y_s^{i,n}-Y_s^{i,p}|(Y^{i,p}_s-Y^{l,p}_s+g_{il})^-ds\biggr]\\ & \qquad+2\E\biggl[\underset{t\geq0}{\sup}\,\bigg|\int_{t}^{+\infty}(Y_s^{i,n}-Y_s^{i,p})(Z_s^{i,n}-Z_s^{i,p})dB_s\bigg|\biggr].
\end{aligned}
\end{equation}
Then, by applying the Burkholder-Davis-Gundy's inequality to the last term of the right hand side of
inequality $(\ref{BDG}),$ we obtain
\begin{equation}\label{BDG2}
\begin{aligned}
&\E\biggl[\underset{t\geq0}{\sup}\,\bigg|\int_{t}^{+\infty}(Y_s^{i,n}-Y_s^{i,p})(Z_s^{i,n}-Z_s^{i,p})dB_s\bigg|\biggr]\\ & \leq C\E\biggl[\int_{t}^{+\infty}|(Y_s^{i,n}-Y_s^{i,p})(Z_s^{i,n}-Z_s^{i,p})|^2dB_s\biggr]^\frac{1}{2},\\ & \leq \frac{1}{2}\E\biggl[\underset{t\geq0}{\sup}\,|Y_t^{i,n}-Y_t^{i,p}|^2\biggr]+C\mathbb{E}\biggr[\int_{t}^{+\infty}|Z_s^{i,n}-Z_s^{i,p}|^2ds\biggr].
\end{aligned}
\end{equation}
Combining $(\ref{BDG})$ and $(\ref{BDG2})$ and taking into consideration $(\ref{lebesgue})$ and $(\ref{Znp})$, we get
\begin{equation}\label{Ynp}
\underset{n,\, p \to +\infty}{lim}\E\biggl[\underset{t\geq0}{\sup}\,|Y_t^{i,n}-Y_t^{i,p}|^2\biggr]=0,
\end{equation}
which means that $(Y^{i,n})_{n\geq1}$ is a Cauchy sequence in $\mathcal{S}^2$. Consequently $Y_t^{i}$ is continuous.\\
\par Now we define $K^{i,n}$ as follows:
\begin{equation}\label{Kin}
K_t^{i,n}:= n\int_{0}^{t}\sum_{l=1}^{m}e^{-rs}(Y^{i,n}_s-Y^{l,n}_s+g_{il})^-ds,\,\,\, \forall t\in[0,+\infty)\mbox{ and }i\in\mathcal{I}.
\end{equation}
From the penalized BSDE $(\ref{pbsde})$, we have
\begin{equation*}
K_t^{i,n}=e^{-r t}Y^{i,n}_t-Y^{i,n}_0+\int_{0}^{t} e^{-rs}f_i(s,X_s,Y_s^{i,n},Z^{i,n}_s) \, \mathrm{d}{s}-  \int_{0}^{t}e^{-rs}Z^{i,n}_s \, \mathrm{d}{B_s},
\end{equation*}
we set
\begin{equation*}
K_t^{i}=e^{-r t}Y^{i}_t-Y^{i}_0+\int_{0}^{t} e^{-r s}f_i(s,X_s,Y^{i}_s,Z^{i}_s) \, \mathrm{d}{s}-  \int_{0}^{t}e^{-rs}Z^{i}_s \, \mathrm{d}{B_s}.
\end{equation*}
Then we deduce immediately by $(\ref{Znp})$ and $(\ref{Ynp})$ that $(K^{i,n})_{n\geq1}$ converges to $K^i$ in $\mathcal{S}^{2}$.
So $K^{i}$ is an increasing process, moreover it is continuous, and so $K^{i}\in\mathcal{K}^2$.\\
Therefore, $\forall i\in\mathcal{I},$ $(Y^i,Z^i,K^i)$ satisfies the first relation in RBSDE $(\ref{rbsde})$. Finally, by the a priori estimate $(\ref{n2})$, we have
\begin{equation*}
\E\biggl[\int_{0}^{+\infty}\big[(Y_s^{i,n}-Y_s^{j,n}+g_{ij})^-\big]^2ds\biggr]\leq \frac{C}{n^2},\; \forall i,j\in\mathcal{I}.
\end{equation*}
Letting $n\longrightarrow +\infty$, we deduce
\begin{equation*}
\E\biggl[\int_{0}^{+\infty}\big[(Y_s^{i}-Y_s^{j}+g_{ij})^-\big]^2ds\biggr]=0,\; \forall i,j\in\mathcal{I}.
\end{equation*}
Hence,
\begin{equation}\label{bar}
Y_s^{i}\geq Y_s^{j}-g_{ij},\; \forall i,j\in\mathcal{I}.
\end{equation}
\underline{\bf{Step 4}}: The minimal boundary condition.\\
At last we need to prove that $\int_{0}^{+\infty} e^{-r s}[Y^{i}_s - \underset{j\in\mathcal{I}^{-i}}{\max}(Y^{j}_s-g_{ij})] \, \mathrm{d}{K^{i}_s} = 0,\;  \forall i\in\mathcal{I}.$ To do so, we should show before that $\int_{0}^{+\infty} e^{-r s}[Y^{i,n}_s - \underset{j\in\mathcal{I}^{-i}}{\max}(Y^{j,n}_s-g_{ij})]^+\, \mathrm{d}{K^{i,n}_s} = 0,\;  \forall i\in\mathcal{I}.$ First we remark that since $K^{i,n}$ is increasing then
\begin{equation}\label{sup}
\int_{0}^{+\infty} e^{-r s}[Y^{i,n}_s - \underset{j\in\mathcal{I}^{-i}}{\max}(Y^{j,n}_s-g_{ij})]^+ \, \mathrm{d}{K^{i,n}_s} \geq 0.
\end{equation}
Actually we have from $(\ref{Kin})$ that, for $i\in\mathcal{I}$
\begin{equation*}
\begin{aligned}
&\int_{0}^{+\infty} e^{-r s}[Y^{i,n}_s - \underset{j\in\mathcal{I}^{-i}}{\max}(Y^{j,n}_s-g_{ij})]^+ \, \mathrm{d}{K^{i,n}_s}\\ & =n\sum_{l=1}^{m}\int_{0}^{+\infty}e^{-r s}[Y^{i,n}_s - \underset{j\in\mathcal{I}^{-i}}{\max}(Y^{j,n}_s-g_{ij})]^+[Y^{i,n}_s -Y^{l,n}_s+g_{il}]^-ds.
\end{aligned}
\end{equation*}
Next in fact it is clear that for $l = i$
\begin{equation*}
[Y^{i,n}_s - Y^{j,n}_s+g_{ij}]^+[Y^{i,n}_s -Y^{l,n}_s+g_{il}]^-=0.
\end{equation*}
On the other hand for $j\neq i,\, l\neq i$ we have
\begin{equation*}
\begin{aligned}
&[Y^{i,n}_s - \underset{j\in\mathcal{I}^{-i}}{\max}(Y^{j,n}_s-g_{ij})]^+[Y^{i,n}_s -Y^{l,n}_s+g_{il}]^-\\ & \leq [Y^{i,n}_s - Y^{l,n}_s+g_{il}]^+[Y^{i,n}_s -Y^{l,n}_s+g_{il}]^-=0.
\end{aligned}
\end{equation*}
Therefore, we deduce that
\begin{equation}\label{inf}
\int_{0}^{+\infty} e^{-r s}[Y^{i,n}_s - \underset{j\in\mathcal{I}^{-i}}{\max}(Y^{j,n}_s-g_{ij})]^+ \, \mathrm{d}{K^{i,n}_s} \leq 0.
\end{equation}
From $(\ref{sup})$ and $(\ref{inf})$ we obtain that
\begin{equation*}
\int_{0}^{+\infty} e^{-r s}[Y^{i,n}_s - \underset{j\in\mathcal{I}^{-i}}{\max}(Y^{j,n}_s-g_{ij})]^+ \, \mathrm{d}{K^{i,n}_s} =0.
\end{equation*}
Now, by applying Lemma $5.8$ in \cite{GP}, the sequence $(\int_{0}^{+\infty} e^{-r s}[Y^{i,n}_s - \underset{j\in\mathcal{I}^{-i}}{\max}(Y^{j,n}_s-g_{ij})]^+ \, \mathrm{d}{K^{i,n}_s})_{n\geq1}$ converges to $\int_{0}^{+\infty} e^{-r s}[Y^{i}_s - \underset{j\in\mathcal{I}^{-i}}{\max}(Y^{j}_s-g_{ij})]^+ \, \mathrm{d}{K^{i}_s}.$ Hence,
\begin{equation*}
\int_{0}^{+\infty} e^{-r s}[Y^{i}_s - \underset{j\in\mathcal{I}^{-i}}{\max}(Y^{j}_s-g_{ij})]^+ \, \mathrm{d}{K^{i}_s} =0.
\end{equation*}
Finally, from $(\ref{bar})$, we conclude that
\begin{equation}\label{lc}
\int_{0}^{+\infty} e^{-r s}[Y^{i}_s - \underset{j\in\mathcal{I}^{-i}}{\max}(Y^{j}_s-g_{ij})] \, \mathrm{d}{K^{i}_s} =0.
\end{equation}
The proof of Proposition $(\ref{prop1})$ is now complete.
\end{proof}
\subsection{Uniqueness}
In this subsection, we prove the uniqueness of $(\ref{keybsde})$ by a verification method. As in Theorem $3.1$ in \cite{HT}, we give a switching  representation property for the solution $Y^{i}$ of $(\ref{keybsde})$, which represents the relationship between this solution and the optimal switching problem.\\
\par In order to state this representation result, first we introduce some notations.\\
Let $a:=(\tau_{n},\zeta_{n})_{n\geq0}$ be an admissible strategy of switching, i.e.,
\begin{itemize}
	\item[-] $(\tau_{n})_{n\geq0}$ is an increasing sequence of stopping times such that $\P(\tau_n<+\infty, \forall n\geq0)=0.$
	\item[-]  $\forall n\geq0,\; \zeta_{n}$ is a random variable with values in $\I$ and $\F_{\tau_n}$-measurable.
	\item[-]  If $(A_t^{a})_{t\geq0}$ is the non-decreasing, $\F$-adapted and càdlàg process defined by
	\begin{equation}
	\forall t\in[0,+\infty),\quad A_t^{a}:=\underset{n\geq1}{\sum}e^{-r\tau_{n}}g_{\zeta_{n-1}\zeta_{n}}1_{[\tau_{n}\leq t]} \mbox{ and } A_{+\infty}^{a}=\underset{t\to+\infty}{\lim}A_t^{a},\; \P-a.s.,
	\end{equation}
	then $\E[(A_{+\infty}^{a})^2]<+\infty$. The quantity $A_{+\infty}^{a}$ stands for the switching cost at infinity when the strategy $a$ is implemented.
\end{itemize}
Next, with an admissible strategy  $a:=(\tau_{n},\zeta_{n})_{n\geq0}$ we associate a state process $(a_t)_{t\geq0}$ defined by
 \begin{equation}
a_t:=\zeta_{0}1_{\{\tau_{0}\}}+\underset{n\geq1}{\sum}\zeta_{n}1_{]\tau_{n-1},\tau_{n}]}, \qquad \forall t\in[0,+\infty).
 \end{equation}
 Finally, for $(i,t)\in\I\times[0,+\infty),$ we also define $\A_{t}^{i}$ the subset of admissible strategies restricted to start in state $i$ at time $t$.\\
 Now, for any $a:=(\tau_{n},\zeta_{n})_{n\geq0}$ which belongs to $\A_{t}^{i}$, let us define the pair of processes
 $(U^a,V^a)$ which belongs to $\mathcal{S}^2\times\mathcal{M}^2$ and which solves the following switched BSDE:
 \begin{equation}\label{switched_bsde}
 e^{-rt}U^{a}_t = \int_{t}^{+\infty} e^{-rs}f_{a_s}(s,X_s,U^{a}_s,V^{a}_s) \, ds -( A^{a}_{\infty} - A^{a}_t)-  \int_{t}^{+\infty}e^{-rs}V^{a}_s \, dB_s,\quad t\in[0,+\infty).
 \end{equation}
 Actually, in setting up $\tilde{U}^a:= e^{-r.}U^{a}-A^{a}$ and $\tilde{V}^a:= e^{-r.}V^{a}$ we remark that BSDE $(\ref{switched_bsde})$ is equivalent to the
 following one:
 \begin{equation}\label{switched_bsde_eq}
 \tilde{U}^a_t = -A^{a}_{\infty}+\int_{t}^{+\infty} e^{-rs}\tilde {f}_{a_s}(s,X_s,\tilde{U}^a_s,\tilde{V}^a_s) \, ds -\int_{t}^{+\infty}\tilde{V}^a_s \, dB_s,\quad t\in[0,+\infty),
 \end{equation}
 where the driver $\tilde{f}_{a_s}$ given by $$\tilde {f}_{a_s}(s,X_s,\tilde{U}^a_s,\tilde{V}^a_s):=f_{a_s}(s,X_s,e^{rs}(\tilde{U}^a_s+A^{a}_s),e^{rs}\tilde{V}^a_s).$$
 Now, since $a$ is admissible and then $\E[(A_{+\infty}^{a})^2]<+\infty$. Therefore, from the result of Chen \cite{chen}, the solution of BSDE $(\ref{switched_bsde_eq})$ exists and is unique.\\
 Hence we deduce that BSDE $(\ref{switched_bsde})$ has a solution in $\mathcal{S}^2\times\mathcal{M}^2$ denoted by $(U^a,V^a)$. \\ \par
  The assumptions required for the uniqueness will be slightly stronger than those needed for existence. We keep the same assumption on $f_i$, and we assume the following
 for the switching costs\\
 \begin{itemize}
 	\item[$\mathbf{[H2']}$] For any $(i,j)\in\I^2$, $g_{ij}$ satisfies the following:
 	\begin{itemize}
 \item[$(i)$]  $g_{ii}\geq0;$
 \item[$(ii)$] $g_{ij}>0,$ for $i\neq j$;
 \item[$(iii)$] for any $(i,j,l)\in\mathcal{I}^3$, such that $i\neq j$ and $j\neq l$, we have $$g_{ij}+g_{jl}> g_{il}.$$
    \end{itemize}
\end{itemize}
 Furthermore,  as in Theorem $3.1$ in \cite{HT}, we give in the Proposition below a switching representation property for the solution $Y^i$ of $(\ref{keybsde})$, which represents the relationship between this solution and the optimal switching
 problem where the aim is to maximize $U^{a}_t$ subject to $a\in\A_{t}^{i}$. The following proposition implies the
 uniqueness of the solution to RBSDE $(\ref{keybsde})$.
\begin{proposition}\label{prop32}
Under $\mathbf{[H1]}$ and $\mathbf{[H2']}$, there exists $a^*\in\A_{t}^{i}$ such that
$$Y^i_t=U^{a^*}_t=\underset{a\in\A_{t}^{i}}{ess \sup}\,U^{a}_t,\quad \forall(i,t)\in\I\times[0,+\infty).$$	
\end{proposition}
\begin{proof}
The proof of Proposition $(\ref{prop32})$ is omitted, since it follows from the same reasoning as in the proof of Theorem $3.1$ in \cite{HT}, where the only difference is that in our framework the horizon is infinite unlike to \cite{HT}.	
\end{proof}
\section{The main result}
Now we give the main result of this paper.
\begin{theorem}\label{theorem41}
	Assume that $\mathbf{[H1]}$ and $\mathbf{[H2']}$  are satisfied. Then the reflected multi-dimensional BSDE $(\ref{rbsde})$  has a unique solution $(Y^i,Z^i,K^i)_{i\in\I}$.
\end{theorem}
\begin{proof}
	 We suppose that, for $i\in\I$, the $i$-th component of the random function $f$ depends on $\vec{y}$.\\
	 Next, we fix $\vec{\Gamma}:=(\Gamma^1,...,\Gamma^m)$ in $[\mathcal{S}^2]^m$ and introduce the operator $\phi:[\mathcal{S}^2]^m\mapsto[\mathcal{S}^2]^m, \vec{\Gamma}\mapsto\vec{Y}:=\phi(\vec{\Gamma})$, where $(\vec{Y},\vec{Z},\vec{K}):=(Y^i,Z^i,K^i)_{i\in\I}\in[\mathcal{S}^2\times\mathcal{M}^2\times\mathcal{K}^2]^m$ is the solution to the following RBSDE, $\forall i\in\mathcal{I}$ and $t \in [0,+\infty)$,
	 \begin{equation}\label{RBSDE41}
	 \begin{cases}
	 e^{-rt}Y^{i}_t = \int_{t}^{+\infty} e^{-rs}f_i(s,X_s,\vec{\Gamma}_s,Z^{i}_s) \, \mathrm{d}{s} + K^{i}_{+\infty} - K^{i}_t-  \int_{t}^{+\infty} e^{-rs}Z^{i}_s \, \mathrm{d}{B_s};\\
	 \lim_{t \to+\infty} e^{-rt}Y^{i}_t=0,\\
	 \forall t \geq 0, \quad e^{-rt}Y^{i}_t \geq e^{-rt} \underset{j\in\mathcal{I}^{-i}}{\max}(Y^{j}_t-g_{ij}),\\
	 \int_{0}^{+\infty} e^{-rs}(Y^{i}_s - \underset{j\in\mathcal{I}^{-i}}{\max}(Y^{j}_s-g_{ij})) \, \mathrm{d}{K^{i}_s} = 0,
	 \end{cases}
	 \end{equation}
	 which exists and is unique thanks to Proposition $(\ref{prop1})$ and $(\ref{prop32})$. Then $\phi$ is well defined and is obviously valued in $[\mathcal{S}^2]^m$.
	 \par Now our objective is to show that $\phi$ is a contraction on $[\mathcal{S}^2]^m$ when endowed with an appropriate equivalent norm.
	 \begin{proposition}\label{prop41}
	 	The operator $\phi$ is a contraction on the Banach space $[\mathcal{S}^2]^m$ endowed with the norm $\|.\|_{2,r}$ defined by:
	 	$$ \|Y\|_{2,r}:=\mathbb{E}\left[\underset{t\geq0}{\sup}\,e^{-rt}|Y_t|^2\right]^{\frac{1}{2}}.$$
	 \end{proposition}
 \begin{proof}
 	 In the same spirit of the proof of Proposition $3.3$ in \cite{CEK} we provide the proof of Proposition $(\ref{prop41})$ only for the sake of completeness, since the proof remains the same even if in our framework we consider
 	an infinite horizon.\\
 	
 	We consider two processes $\vec{\Gamma},\vec{\dot{\Gamma}}\in[\mathcal{S}^2]^m$ such that $\vec{Y}:=\phi(\vec{\Gamma})$ and $\vec{\dot{Y}}:=\phi(\vec{\dot{\Gamma}})$, where $(Y^i,Z^i,K^i)_{i\in\I}$ (respectively, $(\dot{Y}^i,\dot{Z}^i,\dot{K}^i)_{i\in\I}$) is the solution of the RBSDE $(\ref{RBSDE41})$.\\
 	Next, let us introduce the following RBSDE:
 	\begin{equation}\label{hatRBSDE}
 	\begin{cases}
 	e^{-rt}\hat{Y}^i_t = \int_{t}^{+\infty} e^{-rs}\hat{f_i}(s,X_s,\hat{Z}^i_s) \, \mathrm{d}{s} + \hat{K}^i_{+\infty} - \hat{K}^i_t-  \int_{t}^{+\infty}e^{-rs}\hat{Z}^i_s \, \mathrm{d}{B_s};\\
 	\lim_{t \to+\infty} e^{-rt}\hat{Y}^i_t=0,\\
 	\forall t \geq 0, \quad e^{-rt}\hat{Y}^i_t \geq e^{-rt} \underset{j\in\mathcal{I}^{-i}}{\max}(\hat{Y}^j_t-g_{ij}),\\
 	\int_{0}^{+\infty} e^{-rs}\{\hat{Y}^i_s - \underset{j\in\mathcal{I}^{-i}}{\max}(\hat{Y}^j_s-g_{ij})\} \, \mathrm{d}{\hat{K}^i_s} = 0,
 	\end{cases}
 	\end{equation}
 	where $$\hat{f_i}:(t,x,z)\mapsto f_i(t,x,\vec{\Gamma},z)\vee f_i(t,x,\vec{\dot{\Gamma}},z),\quad \forall(i,x,z)\in\I\times[0,+\infty)\times\R^k\times\R^d.$$
 	Once again, by Proposition $(\ref{prop1})$ and $(\ref{prop32})$,  there exists a unique solution to $(\ref{hatRBSDE})$. From Proposition $(\ref{prop32})$, for $(i,t)\in\I\times[0,+\infty)$ and for any $a\in\A_{t}^i$, $\hat{Y}$, $Y^i$ and $\dot{Y}^i$ have the following switching representation property:
 	\begin{equation}\label{UYswitch}
 	\hat{Y}^i_t=\hat{U}^{\hat{a}}_t:=\underset{a\in\A_{t}^{i}}{ess \sup}\,\hat{U}^{a}_t,\quad Y^i_t=\underset{a\in\A_{t}^{i}}{ess \sup}\,U^{a}_t,\mbox{  and  } \dot{Y}^i_t=\underset{a\in\A_{t}^{i}}{ess \sup}\,\dot{U}^{a}_t;
 	\end{equation}
 	where $\hat{a}\in\A_{t}^i$ and $(\hat{U}^a,\hat{V}^a)$, $(U^a,V^a)$ and $(\dot{U}^a,\dot{V}^a)$ are respectively solutions of the following BSDE:
 	$$e^{-rt}\hat{U}^a_t = \int_{t}^{+\infty} e^{-rs}\hat{f}_{a_s}(s,X_s,\hat{V}^{a}_s) \, ds -( A^{a}_{\infty} - A^{a}_t)-  \int_{t}^{+\infty}e^{-rs}\hat{V}^{a}_s \, dB_s,\quad t\in[0,+\infty),$$
 	$$e^{-rt}U^{a}_t = \int_{t}^{+\infty} e^{-rs}f_{a_s}(s,X_s,\vec{\Gamma}_s,V^{a}_s) \, ds -( A^{a}_{\infty} - A^{a}_t)-  \int_{t}^{+\infty}e^{-rs}V^{a}_s \, dB_s,\quad t\in[0,+\infty),$$ and
 	$$e^{-rt}\dot{U}^{a}_t = \int_{t}^{+\infty} e^{-rs}f_{a_s}(s,X_s,\vec{\dot{\Gamma}}_s,\dot{V}^{a}_s) \, ds -( A^{a}_{\infty} - A^{a}_t)-  \int_{t}^{+\infty}e^{-rs}\dot{V}^{a}_s \, dB_s,\quad t\in[0,+\infty).$$
 	Now since, $\forall(i,t,x,z)\in\I\times[0,+\infty)\times\R^k\times\R^d$ we have that $\hat{f_i}(t,x,z)\geq f_i(t,x,\vec{\Gamma},z)$ and $\hat{f_i}(t,x,z)\geq f_i(t,x,\vec{\dot{\Gamma}},z)$, by the classical comparison theorem of BSDE, for any $a\in\A_{t}^i$
 	\begin{equation*}
 	U^{a}_t\leq\hat{U}^a_t \,\mbox{ and }\, \dot{U}^{a}_t\leq\hat{U}^a_t,\quad \forall t\geq 0.
 	\end{equation*}
 	Next, combining this estimates with the representation $(\ref{UYswitch})$, implies that
 	\begin{equation*}
 	Y^{i}_t\leq\hat{Y}^i_t \,\mbox{ and }\, \dot{Y}^{i}_t\leq\hat{Y}^i_t,\quad \forall t\geq 0.
 	\end{equation*}
 	Since $\hat{a}$  is an admissible strategy for
 	the representation $(\ref{UYswitch})$ of $Y^{i}$ and $\dot{Y}^i$,  we deduce that
 	\begin{equation}
 	U^{\hat{a}}_t \leq Y^{i}_t\leq\hat{U}^{\hat{a}}_t \,\mbox{ and }\, \dot{U}^{\hat{a}}_t\leq \dot{Y}^{i}_t\leq\hat{U}^{\hat{a}}_t,\quad \forall t\geq 0.
 	\end{equation}
 	Then, we obtain that
 	\begin{equation}\label{y-doty}
 	|Y^{i}_t-\dot{Y}^{i}_t|\leq|\hat{U}^{\hat{a}}_t-U^{\hat{a}}_t|+|\hat{U}^{\hat{a}}_t-\dot{U}^{\hat{a}}_t|,\quad \forall t\geq 0.
 	\end{equation}
 	We first control the first term on the right hand side of $(\ref{y-doty})$. Using Itô's formula to $e^{-rt}|\hat{U}^{\hat{a}}_t-U^{\hat{a}}_t|^2$, we get
 	\begin{equation}\label{hatU-U}
 	\begin{aligned}
 	 &e^{-rt}|\hat{U}^{\hat{a}}_t-U^{\hat{a}}_t|^2+\int_{t}^{+\infty}e^{-rs}(2r|\hat{U}^{\hat{a}}_s-U^{\hat{a}}_s|^2+|\hat{V}^{\hat{a}}_s-V^{\hat{a}}_s|^2)ds\\ &= 2\int_{t}^{+\infty}e^{-rs}(\hat{U}^{\hat{a}}_s-U^{\hat{a}}_s)(\hat{f}_{\hat{a}_s}(s,X_s,\hat{V}^{\hat{a}}_s)-f_{\hat{a}_s}(s,X_s,\vec{\Gamma}_s,V^{\hat{a}}_s))ds\\&\qquad-2\int_{t}^{+\infty}e^{-rs}(\hat{U}^{\hat{a}}_s-U^{\hat{a}}_s)(\hat{V}^{\hat{a}}_s-V^{\hat{a}}_s)dB_s.
 	\end{aligned}
 	\end{equation}
 	Making use of inequality $|x\vee y-y|\leq|x-y|$ together with the Lipschitz property of $f_.$, we obtain
 	\begin{equation}\label{hatf-f}
 	|\hat{f}_{\hat{a}_s}(s,X_s,\hat{V}^{\hat{a}}_s)-f_{\hat{a}_s}(s,X_s,\vec{\Gamma}_s,V^{\hat{a}}_s)|\leq u(s)(|\vec{\Gamma}_s-\vec{\dot{\Gamma}}_s|+|\hat{V}^{\hat{a}}_s-V^{\hat{a}}_s|).
 	\end{equation}
 	Then,
 	\begin{equation*}
 	\begin{aligned}
 	 &e^{-rt}|\hat{U}^{\hat{a}}_t-U^{\hat{a}}_t|^2+\int_{t}^{+\infty}e^{-rs}(r|\hat{U}^{\hat{a}}_s-U^{\hat{a}}_s|^2+|\hat{V}^{\hat{a}}_s-V^{\hat{a}}_s|^2)ds\\ &\leq \int_{t}^{+\infty}u(s)e^{-rs}|\hat{U}^{\hat{a}}_s-U^{\hat{a}}_s|(|\vec{\Gamma}_s-\vec{\dot{\Gamma}}_s|+|\hat{V}^{\hat{a}}_s-V^{\hat{a}}_s|)ds\\&\qquad-2\int_{t}^{+\infty}e^{-rs}(\hat{U}^{\hat{a}}_s-U^{\hat{a}}_s)(\hat{V}^{\hat{a}}_s-V^{\hat{a}}_s)dB_s.
 	\end{aligned}
 	\end{equation*}
 	Using the elementary inequality $2xy\leq\frac{x^2}{\epsilon}+\epsilon y^2$ for any $x,y\in\R$ and $\epsilon>0$, we get that
 	\begin{equation}\label{lipsch}
 	\begin{aligned}
 	 &e^{-rt}|\hat{U}^{\hat{a}}_t-U^{\hat{a}}_t|^2+\int_{t}^{+\infty}e^{-rs}(r|\hat{U}^{\hat{a}}_s-U^{\hat{a}}_s|^2+|\hat{V}^{\hat{a}}_s-V^{\hat{a}}_s|^2)ds\\ &\leq \int_{t}^{+\infty}(2u^2(s)+\epsilon^{-1})e^{-rs}|\hat{U}^{\hat{a}}_s-U^{\hat{a}}_s|^2ds+\epsilon\int_{t}^{+\infty}u^2(s)e^{-rs}|\vec{\Gamma}_s-\vec{\dot{\Gamma}}_s|^2ds\\&\qquad+\frac{1}{2}\int_{t}^{+\infty}e^{-rs}|\hat{V}^{\hat{a}}_s-V^{\hat{a}}_s|^2ds-2\int_{t}^{+\infty}e^{-rs}(\hat{U}^{\hat{a}}_s-U^{\hat{a}}_s)(\hat{V}^{\hat{a}}_s-V^{\hat{a}}_s)dB_s.
 	\end{aligned}
 	\end{equation}
 	Choosing $r \geq 2u^2(s)+\epsilon^{-1}$ and taking expectation at $t=0$, we have
 	\begin{equation}\label{V2intgra}
 	\E\biggl[\int_{0}^{+\infty}e^{-rs}|\hat{V}^{\hat{a}}_s-V^{\hat{a}}_s|^2ds\biggr]\leq 2\epsilon\E\biggl[\int_{0}^{+\infty}u^2(s)e^{-rs}|\vec{\Gamma}_s-\vec{\dot{\Gamma}}_s|^2ds\biggr].
 	\end{equation}
 	Going back to $(\ref{lipsch})$ and applying Burkholder-Davis-Gundy's inequality, we obtain
 	\begin{equation*}
 	\begin{aligned}
 	\E\biggl[\underset{t\geq0}{\sup}\,e^{-rt}|\hat{U}^{\hat{a}}_t-U^{\hat{a}}_t|^2\biggr]&\leq 2\epsilon\E\biggl[\int_{0}^{+\infty}u^2(s)e^{-rs}|\vec{\Gamma}_s-\vec{\dot{\Gamma}}_s|^2ds\biggr]+\frac{1}{2}\E\biggl[\underset{t\geq0}{\sup}\,e^{-rt}|\hat{U}^{\hat{a}}_t-U^{\hat{a}}_t|^2\biggr]\\&\qquad+\frac{C^2}{2}\E\biggl[\int_{0}^{+\infty}e^{-rs}|\hat{V}^{\hat{a}}_s-V^{\hat{a}}_s|^2ds\biggr].
 	\end{aligned}
 	\end{equation*}
 	Taking into account inequality $(\ref{V2intgra})$, we get
 	\begin{equation}
 		\E\biggl[\underset{t\geq0}{\sup}\,e^{-rt}|\hat{U}^{\hat{a}}_t-U^{\hat{a}}_t|^2\biggr]\leq C_{\epsilon}\E\biggl[\int_{0}^{+\infty}u^2(s)e^{-rs}|\vec{\Gamma}_s-\vec{\dot{\Gamma}}_s|^2ds\biggr],
    \end{equation}
     where $C_{\epsilon}$ denotes a constant which depends on $\epsilon$  and may vary from line to line.\\
 	Following the same method, we can get a similar estimate for $\E[\underset{t\geq0}{\sup}\,e^{-rt}|\hat{U}^{\hat{a}}_t-\dot{U}^{\hat{a}}_t|^2]$.\\
 	Next going back to $(\ref{y-doty})$ and taking in consideration the estimates obtained above, we deduce
 	\begin{equation}
 	\E\biggl[\underset{t\geq0}{\sup}\,e^{-rt}|Y^{i}_t-\dot{Y}^{i}_t|^2\biggr]\leq C_{\epsilon}\E\biggl[\int_{0}^{+\infty}u^2(s)e^{-rs}|\vec{\Gamma}_s-\vec{\dot{\Gamma}}_s|^2ds\biggr].
 	\end{equation}
 	Since the last inequality holds true for an $(i,t)\in\I\times[0,+\infty)$, we have
 	$$\|\phi(\vec{\Gamma})-\phi(\vec{\dot{\Gamma}})\|^2_{2,r}\leq C_{\epsilon}\int_{0}^{+\infty}u^2(s)ds\|\vec{\Gamma}-\vec{\dot{\Gamma}}\|^2_{2,r}$$
 	From $\int_{0}^{+\infty}u^2(s)ds<C_{\epsilon}^{-1}$, $\phi$ is a contraction on $[\mathcal{S}^2]^m$, which concludes the proof of the proposition.
 \end{proof}
As a consequence, there exists a unique fixed point in $[\mathcal{S}^2]^m$ for $\phi$, which is the unique solution of RBSDE $\eqref{rbsde}$.
 \end{proof}

\end{document}